\documentclass[leqno,10pt, twocolumn]{article}
\usepackage{amsfonts}
\usepackage{latexsym}
\usepackage{amsmath}
\usepackage{epsfig}
\usepackage{ascmac}
\usepackage{amssymb}

\newtheorem{thm}{Theorem}[section]
\newtheorem{lem}[thm]{Lemma}
\newtheorem{cor}[thm]{Corollary}

\def\0{{\bf 0}}

\def\a{{\bf a}}

\def\x{{\mib x}}
\def\y{{\mib y}}

\def\D{\mathbb{D}}
\def\N{\mathbb{N}}

\def\R{\mathbb{R}}

\def\X{{\bf X}}
\def\Y{{\bf Y}}

\def\cD{\mathcal{D}}
\def\cE{\mathcal{E}}

\def\cP{\mathcal{P}}

\def\cX{\mathcal{X}}

\def\M{{\mathfrak M}}

\def\Ai{{\rm Ai}}

\newcommand{\mib}[1]{\mbox{\boldmath $#1$}}

\newcommand{\qed}{\hbox{\rule[-2pt]{3pt}{6pt}}}

\title{{\bf  Cores of Dirichlet forms %
related to \\Random Matrix Theory} }
\author{ Hirofumi Osada
\footnote{
Faculty of Mathematics, Kyushu University, Fukuoka, 819-0395, JAPAN, email:osada@math.kyushu-u.ac.jp
}
, \quad Hideki Tanemura
\footnote{
Department of Mathematics and Informatics, Faculty of Science, Chiba University, 
1-33 Yayoi-cho, Inage-ku, Chiba 263-8522, JAPAN, 
email:tanemura@math.s.chiba-u.ac.jp
}
}

\begin{document}
\maketitle
\begin{abstract}
We prove the sets of polynomials on configuration spaces are cores of Dirichlet forms describing interacting Brownian motion in infinite dimensions. Typical examples of these stochastic dynamics are Dyson's Brownian motion and Airy interacting Brownian motion. 
Both particle systems have logarithmic interaction potentials, and naturally arise from random matrix theory. The results of the present paper will be used in a forth coming paper to prove the identity of the infinite-dimensional stochastic dynamics related to the random matrix theories constructed by apparently different methods: the method of space-time correlation functions and that of stochastic analysis.   
\end{abstract}


\section{Introduction.}

In random matrix theory,  one of the main issues is to clarify 
the distribution of the eigenvalues and 
its asymptotic behavior as the size of the matrices goes to infinity. 
The prototypes of random matrices are Gaussian ensembles, 
 divided into three classes and 
called Gaussian orthogonal/unitary/symplectic ensembles (GOE/GUE/GSE), 
according to their invariance under conjugation by orthogonal/unitary/symplectic groups. 

The eigenvalue distributions of Gaussian random matrices 
of $ N\times N $ size  are then given as 
\begin{equation}
\nonumber
\check{m}_\beta^N(d\x_N)=\frac{1}{Z}h_N(\x_N)^\beta 
e^{-\frac{\beta}{4}|\x_N|^2}
d\x_N,
\end{equation}
where $d\x_n=dx_1dx_2\cdots dx_N$, $\x_N=(x_i)\in\mathbb{R}^N$, and 
$$\displaystyle{h_N(\x_N)=\prod_{i<j}^N|x_i-x_j|}
.$$  
Here and after $Z$ denotes the normalizing constant. 
The GOE, GUE, and GSE correspond to inverse temperature 
$\beta=1,2$ and $4$, respectively \cite{Meh04, AGZ10}.

The celebrated Wigner's theorem asserts that the empirical measure 
\begin{align}\label{:0a}&
\displaystyle{\frac{1}{N}\sum_{j=1}^N \delta_{{x_j}/{\sqrt{N}}}}
.\end{align}
of the eigenvalues under the distribution $\check{m}_\beta^N(d\x_N)$ converges to the semicircle law 
$ \varsigma (x) dx$ as $N\to\infty$ \cite{AGZ10, Meh04}, where 
$ \varsigma (x) dx $ is the probability on $ \mathbb{R}$ 
such that 
\begin{align}\label{:0}&
\varsigma (x) dx = \frac{1}{2\pi}\sqrt{4-x^2}{\bf 1}_{[-2,2]}(x)dx
.\end{align}

For a countable subset $ \{x_n\}$, we call the $ \sigma $-finite measure 
$ \xi = \sum_n \delta_{x_n}$ a configuration if it becomes a Radon measure. 
The set of all configurations on $ \mathbb{R}$ is a Polish space equipped with the vague topology, and is called the configuration space over $ \mathbb{R}$. 
We call $ \xi $ unlabeled particles, and $ {\mib x}=(x_1,x_2,\ldots )$ labeled particles. 

Note that for $ N$-particle systems, there exists an obvious bijection between the distribution of the unlabeled $ N$-particles and the {\em symmetric}  distribution of the labeled $ N$-particles, where $ N\in\N $. 
We note that this is not the case for infinite particle systems.

For a given distribution $ \mu $ of $ N$-unlabeled particles, 
we denote by $ \check{\mu}$ the symmetric density of the associated 
$ N$-labeled particles in the sequel. 

%

To examine the behavior of the distribution of the configuration 
$ \xi = \sum_{i=1}^{N} \delta_{x_j}$ under 
$\check{m}_\beta^N(d\x_N)$, $N\to\infty$ there are two typical scalings, 
 called the {\it bulk} and the {\it soft edge}. 
The former corresponds to the scaling such that $y_j= \sqrt{N}x_j$, and the distribution of $\{y_j\}_{j=1}^N$ under $\check{m}_\beta^N(d\x_N)$ 
 is given by 
\begin{equation}
\nonumber
 \check{\mu}_{{\sf bulk}, \beta}^N(d\y_N)=\frac{1}{Z}h_N(\y_N)^\beta 
e^{-\frac{\beta}{4N}|\y_N|^2}
d\y_N.
\end{equation}
The latter corresponds to the scaling such that $y_j= N^{1/6}(x_j-2\sqrt{N})$, and the distribution of 
$\{y_j\}_{j=1}^N$ under $\check{m}_\beta^N$ is given by 
\begin{equation}
\nonumber
\check{\mu}_{{\sf soft}, \beta}^N(d\y_N)=\frac{1}{Z}
h_N(\y_N)^\beta
e^{\frac{-\beta}{4N^{1/3}}|\y_N-2N^{2/3}{\bf 1}_N|^2}
d\y_N,
\end{equation}
where ${\bf 1}_N=(1,1,\dots,1)\in\R^N$.

Let $\beta=2$. Then the limit of $\mu_{{\sf bulk},2}^N$ 
 is the determinantal random point field $\mu_{\sin, 2}$ with sine kernel 
\begin{equation}
\label{def:sin}
K_{\sin, 2}(x,y)=\frac{\sin ( x-y)}{\pi(x-y)} 
,\end{equation}
and the limit of $\mu_{{\sf soft},2}^N$ 
is the determinantal random point field $\mu_{\Ai, 2}$ with Airy kernel
\begin{equation}
\label{def:Ai}
K_{\Ai, 2}(x,y)=\frac{\Ai(x)\Ai'(y)-\Ai'(x)\Ai(y)}{x-y},
\end{equation}
where $\Ai$ denotes the Airy function and $\Ai'$ its derivative \cite{ST03, Sos00, Meh04}.
It is proved that these random point fields are 
quasi-Gibbsian in \cite{o.rm, o.rm2}.

We consider the dynamical scaling limit corresponding to the static limit mentioned above. 
For this we introduce the associated stochastic dynamics describing the time evolution of $ N$-particle systems.

Let $\X^N(t)=(X_j^N(t))_{j=1}^{N}$ be the solution of the SDE
\begin{align}&
\label{:1}
&&dX_j^N(t)=dB_j(t)+\sum_{k=1,k\not=j}^N \frac{dt}{X_j^N(t)-X_k^N(t)}
\end{align}
or the SDE with Ornstein-Uhlenbeck's type drifts 
\begin{align} \label{:2} &
 dX_j^N(t)=dB_j(t) - \frac{1}{2N} X_j^N(t)dt 
\\ \notag & \quad \quad \quad \quad \quad \quad \quad + 
\sum_{k=1,k\not=j}^N \frac{dt}{X_j^N(t)-X_k^N(t)}
.\end{align}
These are called Dyson's Brownian motion model with $\beta=2$, 
or simply the Dyson model \cite{Dys62}. 
The solution of \eqref{:2} is a natural reversible stochastic dynamics with respect to 
$  \check{\mu}_{{\sf bulk},2}^N $, and that of \eqref{:1} is also natural but has no invariant probability measures. 
Both have the same $ N$-limit 
as we see below.

Let $ \Xi _{\sin}^N $ be an unlabeled process defined by 
$$\Xi _{\sin}^N (t)=\sum_{j=1}^N \delta_{X_j^N(t)}. $$
Suppose that the distribution of $\Xi _{\sin}^N (0)$ is 
$ \mu_{{\sf bulk},2}^N $. 
Then $\Xi _{\sin}^N $ converges in distribution to 
the process $\Xi_{\sin} $ whose generating function 
$$
{\Psi}_{\sin}^{\bf t}[{\bf f}]
\equiv {\bf E} \left[\exp \left\{ \sum_{m=1}^{M} 
\int_{\R} f_m(x) \Xi_{\sin}(t_m, dx) \right\} \right],
$$
$0\le t_1 \le t_2 \le \cdots \le t_M$, $f_m \in C_0(\R)$, $1\le m\le M$, is represented by the Fredholm determinant 
\begin{equation}
 \mathop{\rm Det}_
{\substack
{(s,t)\in \{t_1, \dots, t_M \}^2, \\
(x,y)\in \R^2}
}
\Big[\delta_{s t} \delta(x-y)
+ {\bf K}_{\sin}(s, x; t, y) \chi_{t}(y) \Big],
\nonumber
\end{equation}
with $\chi_{t_m}= e^{f_m}-1$ and 
the {\it extended sine kernel} ${\bf K}_{\sin}$ \cite{Spo87, KT07}:
\begin{align*}
&{\bf K}_{\sin}(s,x;t,y) 
\\
&= 
\begin{cases}
\displaystyle{
\frac{1}{\pi}\int_{0}^{1} du \, e^{ u^2 (t-s)/2} 
\cos \{ u (y-x)\} 
},
& t\ge s
\cr
\displaystyle{
- \frac{1}{\pi}\int_{1}^{\infty} du \, 
e^{ u^2 (t-s)/2} \cos \{ u (y-x) \} 
},
& t<s.
\end{cases} 
\end{align*}

For the soft edge scaling, we suppose that the distribution of 
$\X^{N}(0)$ is $\check{\mu}_{{\sf bulk},2}^N$, 
and introduce the process $\Y^N $ defined by 
\begin{equation}
\label{SOft_edge_scaling}
\Y^N(t)= \frac{1}{N^{1/3}}\X^N(N^{2/3}t) - 2N^{2/3}-N^{1/3}t+ \frac{t^2}{4}
\end{equation} 
corresponding to \eqref{:1}, and 
\begin{align}\label{:7a}& \!\!\!\!\!\!
\Y^N(t)=\frac{1}{N^{1/3}}\X^N(N^{2/3}t) - 2N^{2/3} 
\end{align}
corresponding to \eqref{:2}. 
Then the unlabeled processes 
$$ \Xi _{\mathrm{Ai}}^N (t) = \sum_{j=1}^N \delta_{Y_j^N(t)}$$
converge in distribution to the process $\Xi_{\Ai} $ 
whose generating function is represented by the Fredholm determinant with 
the {\it extended Airy kernel} ${\bf K}_{\Ai}$ \cite{PS02, Joh03, NKT03, KT07}:
\begin{align*}
&{\bf K}_{\Ai}(s,x;t,y) 
\\
&=\begin{cases}
\displaystyle{
\int_{0}^{\infty} du \, e^{-u(t-s)/2} \Ai(u+x) \Ai(u+y),
}
&t \geq s
\cr 
\displaystyle{
- \int_{-\infty}^{0} d u \, e^{-u(t-s)/2} \Ai(u+x) \Ai(u+y),
} 
&t < s.
\end{cases}
\end{align*}
From the fact that 
${\bf K}_{\sin}(s,x; s,y)=K_{\sin}(x,y)$ and ${\bf K}_{\Ai}(s,x; s,y)=K_{\Ai}(x,y)$, we see that the processes $\Xi_{\sin}$ and $\Xi_{\Ai}$ are reversible with respect to $\mu_{\sin, 2}$ and $\mu_{\Ai, 2}$, respectively.  

The scaling limits above are based on the convergence of the associated space-time correlation functions by the determinatal structures \cite{KT07}. 
There exists another method of constructing infinite volume stochastic dynamics based on stochastic analysis. 
In \cite{o.rm, o.rm2},
unlabeled diffusion processes $\widehat{\Xi}_{\sin}$ and $\widehat{\Xi}_{\Ai}$ 
with reversible measures $\mu_{\sin, 2}$ and $\mu_{\Ai, 2}$ 
are constructed through the Dirichlet form technique.

Let $\X(t)=(X_j(t))_{j\in\N}$ be a labeled process associated with 
$\widehat{\Xi}_{\sin}(t)= \sum_{j\in\N}\delta_{X_j(t)}$. In  \cite{o.isde} it is proved that the process $\X =(X_j)_{j\in\N}$
 solves the infinite-dimensional stochastic differential equation (ISDE)
\begin{equation}
\tag{\bf sin}
dX_j(t)=dB_j(t)+\sum_{k=1,k\not=j}^{\infty}\frac{dt}{X_j(t)-X_k(t)}.
\end{equation}
In \cite{ot.1} we prove that a labeled process $\Y (t)=(Y_j(t))_{j\in\N}$ 
associated with 
$\widehat{\Xi}_{\Ai}(t)= \sum_{j\in\N}\delta_{Y_j(t)}$ 
solves the ISDE 
\begin{align}
\tag{\bf Ai} 
&&\quad \quad dY_j(t)=dB_j(t)
\\
&&+\lim_{r\to\infty}\{\sum_{\substack{k\not=j \\ |Y_k(t)|<r}}^{\infty}\frac{1}{Y_j(t)-Y_k(t)}-\int_{-r}^r \frac{\widehat{\rho}(x)dx}{-x}\}dt,
\nonumber
\end{align}
where
$$
\widehat{\rho}(x)=\frac{{\bf 1}_{(-\infty,0)}(x)}{\pi}\sqrt{-x}.
$$

These two approaches are fundamentally different. 
Hence it is significant to prove that 
the resulting stochastic dynamics are the same. 
From the former construction we can obtain {\em quantitative} information of the limit stochastic dynamics through the calculation of space-time correlation functions; while from the latter we deduce many {\em qualitative} properties of the sample paths of the labeled diffusions through the ISDE representation 
of the processes. 

Recently, we have proved the coincidence of these pairs of stochastic dynamics 
$ \Xi_{\sin}$ and $  \widehat{\Xi}_{\sin} $, and also 
$ \Xi_{\Ai}$ and $  \widehat{\Xi}_{\Ai} $,  
through the following three steps: 
Below $ \star $ denotes  sin or Ai for the sake of brevity. 

\vskip 1mm
\noindent (i) $\Xi_{\star}$  has the strong Markov property.
\vskip 1mm
\noindent (ii) The Dirichlet forms associated with 
$\Xi_{\star}$ and $\widehat{\Xi}_{\star}$ are both 
extensions of the closable form $(\cE^{\mu_{\star}}, {\cal P})$. 
Here $ \cE^{\mu_{\star}}$ are given by \eqref{def:E} and 
${\cal P}$ is the set of polynomials functions on $ \M $ 
defined in (\ref{pol}) later. 

\vskip 1mm
\noindent (iii) The labeled process associated with 
$\Xi_{\star}$ and $\widehat{\Xi}_{\star}$ are solutions of 
the ISDE \thetag{$\star$}, and 
the ISDE \thetag{$\star$} has strong uniqueness. 

\medskip

The claim (i) is proved in \cite{t14}, and 
 \thetag{ii} is in this article. 
The claim (iii) is proved in \cite{ot.1,ot.2} 
partly through the result in this paper. 
Putting these together, we will 
complete the proof of 
$\Xi_{\star} = \widehat{\Xi}_{\star}$ for 
$ \star \in \{\mathrm{sin},\mathrm{Ai}\}$ 
in a forthcoming paper. 

\medskip

\section{Preliminaries}
Let $S$ be a closed subset of $\R^d$ such that the interior $S_{int}$ is a connected open set and that its closure $\overline{S_{int}}$ equals $S$.
Let $\M=\M(S)$ be the configuration space over $S$ of unlabeled particles, the set of non-negative integer valued Radon measures on $S$. The space $\M$ is a Polish space endowed with the vague topology. 
An element $\xi$ of $\M$ can be represented as
$ \xi= \sum_{j\in \Lambda}\delta_{x_j}$
for some countable set $\Lambda$, and the restriction of $\xi$ on a subset $A$ of $S$ is denoted by $\xi_A = \xi (\cdot \cap A)$. 
A function $f$ on $\M$ is called {\it local} if $f(\xi)=f(\xi_K)$ for some compact set $K$.

We write $ \xi_K = \sum_{j=1}^k \delta_{x_j}$. 
For a local function $f$ with $f(\xi)=f(\xi_K)$ 
we introduce the functions $\check{f}_k$ 
on $S^k$, $k\in\N_0\equiv \{0\}\cup \N $ defined by 
$\check{f}_0=f(\emptyset)$, where $\emptyset$ is the null configuration,  and  by, for $ k\in \N $, 
$$
\check{f}_k(\x_k)=f\left(\sum_{j=1}^k \delta_{x_j} \right)   
\text{ for }
\x_k \in K^k 
.$$
We extend the domain of $ \check{f}_k(\x_k)$ to $ S^k\backslash K^k$ 
by the consistency coming from $ f (\xi) = f (\xi_K)$. 
Hence $\check{f}_k$, $k\in\N_0$, satisfy the consistency relation
\begin{equation}
\label{consistency}
\check{f}_{k+1}(\x_k, y)=\check{f}_{k}(\x_k), \ \x_k\in S^k, \ y\notin K.
\end{equation}
The infinite sequence given by
\begin{equation}
\label{inf_seq}
(\check{f}_0, \check{f}_1(x_1),\check{f}_2(x_1,x_2),\dots) = (\check{f}_k(\x_k))_{k\ge 0}
\end{equation}
is a representation of the local function $f$. 

A local function $f$ is called {\it smooth} if the $\check{f}_k$ are smooth 
for $k\in\N_0$. 
We denote by ${\cal D}_{\infty}$ the set of all local smooth functions on $\M$.

Set for $\x_k=(x_i)_{i=1}^k \in S^k$, $k\in\N_0$, 
$  f,g \in \cD_\infty $
$$
\D(f,g)(\x_k)= \frac{1}{2}
\sum_{i=1}^k \sum_{j=1}^d 
\frac{\partial \check{f}_k(\x_k)}{\partial x_{ij}}
\frac{\partial \check{g}_k(\x_k)}{\partial x_{ij}}
,$$
where $x_i=(x_{i1}, x_{i2}, \dots, x_{id})$.
For given $f,g\in\cD_\infty$, the right hand side is a permutation invariant function, and the square field $\D(f,g)$ can be regarded as a local function with variable $\xi=\sum_{i\in \N}\delta_{x_i}\in \M$.

For a probability $\mu$ on $\M$, 
$ L^2(\M,\mu)$ denotes the space of square integrable functions on $\M$ with
the inner product $\langle \cdot, \cdot \rangle_\mu$ and the norm $\|\cdot\|_{L^2(\M,\mu)}$. 
We consider the bilinear form $(\cE^\mu, \cD_\infty^\mu)$ on $L^2(\M,\mu)$ defined as 
\begin{eqnarray}
&&\cE^\mu (f,g)= \int_{\M}\D(f,g)d\mu,
\label{def:E}
\\
&&\cD_\infty^\mu=\{f\in \cD_\infty : \|f\|_1^2 <\infty \}
\label{def:D}
,\end{eqnarray}
where $\|f\|_1^2 \equiv \cE^\mu (f,f)+\|f \|_{L^2(\M,\mu)}^2$. 
A function $F$ on $\M$ is called a polynomial function if $F$ is given as
\begin{equation}
\label{pol}
F(\xi) =Q \left( \langle \phi_1,\xi\rangle, \langle \phi_2,\xi\rangle, \dots, \langle \phi_\ell,\xi\rangle\right) 
\end{equation}
with $\phi_k \in C_0^{\infty}(\R^d)$ and 
a polynomial function $Q$ on $\R^\ell$, where 
$\langle \phi, \xi\rangle = \int_{\R^d}\phi(x) \xi (dx) $ and 
$ C_0^{\infty}(\R^d)$ is the set of smooth functions with compact support. 

We denote by  $\cP $ the set of all polynomial functions on $\M$, 
and by $ \cP _0 $ if we replace the set of polynomials $ Q $ on $\R^\ell$ 
by $ C_b^{\infty} (\R^\ell)$, 
the set of bounded smooth functions with bounded derivaives of any order. 
It is obvious that each element of $ \cP $ and $ \cP _0 $ 
 is a local smooth function. 

The closability and the quasi-regularity of the bilinear form 
$(\cE^\mu, \cD_\infty^\mu)$ have been proved 
 in \cite{o.dfa, o.rm, o.rm2}, while those of  $(\cE^\mu, \cP _0 ) $ 
 in \cite{Y96, AKR}. 
 
Let $ \cP ^{\mu }$, $ \cP _0^{\mu }$, and $ {\cal D}^{\mu }$ 
be the closures with respect to $ \|f\|_1 $ of 
$ \cP $, $ \cP _0 $, and $ {\cal D}_{\infty}^{\mu }$, respectively. 
We see that 
 $ \cP \subset {\cal D}_{\infty}^{\mu } $ and 
 $ \cP _0 \subset {\cal D}_{\infty}^{\mu } $ 
under the mild assumption \thetag{A.0} below, and 
hence we obtain that 
\begin{align}\label{:25}&
\cP ^{\mu } \subset {\cal D}^{\mu } , \quad 
 \cP _0^{\mu } \subset {\cal D}^{\mu } 
.\end{align}
Then we deduce from Theorem \ref{th:main} below that   
\begin{align}\label{:26}&
\cP ^{\mu } = \cP _0^{\mu } = {\cal D}^{\mu } 
.\end{align}

The construction of unlabeled diffusions of 
interacting Brownian motion in infinite dimensions 
through the Dirichlet form approach was initiated by \cite{o.dfa}. 
Later \cite{Y96, AKR} also used this approach but with different cores 
$\cP _0 $ under a more restrictive assumptions 
of interaction potentials than \cite{o.dfa}.  
The identity \eqref{:26} above proves that these diffusions are the same. 

We refer to \cite{MR} and \cite{FOT} for the notion of quasi-regularity and Dirichlet forms.

\section{Main results.}

We call a function $\rho^n$ the $n$-correlation function of $\mu$ with respect to the Lebesgue measure if $\rho^n : S^n \to \R$ is a permutation invariant function such that
\begin{eqnarray}
&&\int_{A_1^{k_1}\times\cdots\times A_m^{k_m}}\rho^n(x_1,\dots,x_n)dx_1\cdots dx_n
\nonumber\\
&&=\int_\M\prod_{i=1}^m \frac{\xi(A_i)!}{(\xi(A_i)-k_i)!}d\mu(\xi)
\nonumber
\end{eqnarray}
for any sequence of disjoint bounded subsets $A_1,\dots,A_n \subset S$ and a sequence of natural numbers $\{k_i\}$ with $k_1+\cdots+k_m=n$.
We assume the following conditions on the probability measure $\mu$ on $\M$:
\vskip 3mm
\noindent
\thetag{A.0} \quad The measure $\mu$ has an $n$-correlation function $\rho^n$ for each $n\in\N$ with  $\rho^n\in L^p(S_r^n,d\x_n)$ 
for all $r\in\N$ for some $1<p\le \infty$. 
Here $S_r=\{x\in S : |x|< r\}$.

\vskip 3mm
\noindent
\thetag{A.1} \quad $(\cE^\mu,\cD_\infty^\mu)$ is closable on 
$ L^2(\M,\mu)$ and its closure 
$(\cE^\mu, \cD^\mu)$ is a quasi-regular Dirichlet form.
\vskip 3mm
\noindent

\medskip 
From \thetag{A.0} we easily deduce \eqref{:25}. 
Hence from \thetag{A.1} we see that 
$(\cE^{\mu}, \cP )$ and $(\cE^{\mu}, \cP_0 )$ are closable on $ L^2(\M,\mu)$ 
as well as $(\cE^\mu,\cD_\infty^\mu)$. 
Let $(\cE^{\mu}, \cP^{\mu})$, $ (\cE^{\mu}, \cP_0^{\mu}) $, and 
$ (\cE^{\mu}, \cD^{\mu})$ be their closures as before. 

The main result of this paper is the following.
\vskip 3mm

\begin{thm}\label{th:main}
Suppose that $\mu$ satisfies {\rm (A.0)} and {\rm (A.1)}.
Then 
$(\cE^{\mu}, \cP^{\mu}) =(\cE^{\mu}, \cP_0^{\mu}) 
= (\cE^{\mu}, \cD^{\mu})$.
\end{thm}

\medskip

In \cite{o.rm, o.rm2} the sufficient conditions (A.0) and (A.1) were given: 
 if $\mu$ is a $(\Phi,\Psi)$-quasi-Gibbs measure, with Borel measurable functions $\Phi: S \to \R\cup\{\infty\}$ and 
 $\Psi: S\times S \to \R\cup\{\infty\}$ satisfying 
\begin{eqnarray}
\label{:16}
&&c^{-1} \Phi_0(x) \le \Phi(x) \le c\Phi_0(x)
\\
\label{:17}
&&c^{-1} \Psi_0(x-y) \le \Psi(x,y) \le c\Psi_0(x-y)
\end{eqnarray}
for a positive constant $c$, and 
some upper semi-continuous functions $\Phi_0$ and $\Psi_0$ being 
 locally bounded from below and with compact core 
 $\{x: \Psi_0(x)=\infty\}$, then 
(A.0) and (A.1) are satisfied. 

\medskip 

We next apply Theorem \ref{th:main} to the stochastic dynamics 
arising from the random matrix theory. 

Since the processes $\Xi_{\sin}$ and $\Xi_{\Ai}$ have the strong Markov property \cite{t14}, they are associated with the quasi-regular Dirichlet forms $({\cal E}_{\sin}, \cal{D}_{\sin})$ and $({\cal E}_{\Ai}, \cal{D}_{\Ai})$, respectively.
Thus we have the desired result as a corollary of Theorem \ref{th:main}. 

\medskip

\begin{cor} \label{c:3}
Let $\star \in \{ \sin , \Ai\}$. 
The Dirichlet forms associated with $\widehat{\Xi}_{\star}$ and 
 $\Xi_{\star}$ are both extensions of the Dirichlet form 
 $(\cE^{\mu_{\star}}, {\cal P}^{\mu_{\star}})$. 
Furthermore, 
$ \cP $ is a core of the Dirichlet form 
$(\cE^{\mu_{\star}}, \cD^{\mu_{\star}})$, and 
\begin{align*}&
\cD^{\mu_{\star}} \subset {\cD}_{\star }
.\end{align*}
\end{cor}

\medskip

\noindent {\em Proof}. 
It is proved in Proposition 7.2 of \cite{KT07} that 
$${\cal E}^{\mu_{\star }} (f,g) = 
{\cal E}_{\star } (f,g),
\quad f,g \in {\cal P}
.$$
From this we deduce the first claim. 
The second is a direct consequence of Theorem \ref{th:main}. 
\qed

\medskip

\section{Proof of Theorem \ref{th:main}.}
For simplicity we only prove the case $ S = \R ^d$. 
We set 
$${\bf A}=\{\a=\{a_r\}_{r\in\N} : a_r\in\N, a_r\le a_{r+1}, r\in\N\}
.$$
For $\a=\{a_r\}\in {\bf A}$, let 
$$
\M[\a]=\{\xi \in \M : \xi (S_r)\le a_r, \mbox{ for all } r\}. 
$$
Then $\M[\a]$ is compact in $ \M $ endowed with the vague topology.
We introduce a cut off function ${\cal X}$ of $\M [\a]$ as follows
\begin{eqnarray}
\cX[\a](\xi)&=& h \circ d_{\a}(\xi),
\nonumber
\\
d_{\a}(\xi)&=&\left\{ 
\sum_{r=1}^\infty \sum_{j\in J_{r,\xi}}(r-|x_j(\xi)|)^2 
\right\}^{1/2}.
\nonumber
\end{eqnarray}
Here $\{x_j(\xi)\}$ is a sequence in $ S $ such that $\xi=\sum\delta_{x_j(\xi)}$, $|x_j(\xi)|\le |x_{j+1}(\xi)|$ for all $j$,
and 
$$
J_{r,\xi}=\{j: j>a_r, \ x_j(\xi )\in S_r\}.
$$
Furthermore, $ h : \R \to [0,1]$ is a function defined by 
\begin{equation}
 h (t)  = 
\begin{cases}
1, &t\in (-\infty,0),
\\
1-t, &t\in [0,1],
\\
0, & t\in (1,\infty).
\end{cases}
\end{equation}
Note that $d_{\a}(\xi)=0$ and ${\cal X}_{\a}(\xi)=1$ 
if $J_{r,\xi}=\emptyset$ for all $r\in\N $. 
The following is proved in \cite[Lemma 2.5]{o.dfa}.

\begin{lem}
\label{l:41}
For any $f\in\cD_\infty^\mu$ and $\varepsilon >0$
we can take $\a\in {\bf A}$ such that $\cX[\a]f\in\cD$ and 
that 
$$
\| (1-\cX[\a]) f\|_1 <\varepsilon.
$$
\end{lem}

Let $\psi$ be a smooth function on $\R$ with support in $[-1,1]$ such that
$\int_\R \ \psi (x) dx=1$.
Then we put $\psi_N(x)=N\psi(Nx)$.
For $ g \in C_0^{\infty} (\R^L) $ with support in $[-r,r]^L$ 
we associate the following function
\begin{eqnarray}
\nonumber
g_N(\x_L)&=&  \sum_{j_1=1}^{N}\cdots\sum_{j_L=1}^N 
g\left(\frac{2j_1 r}{N} - r ,\dots,\frac{2j_L r}{N}- r \right)
\nonumber\\
&&\qquad\qquad\times
\prod_{\ell=1}^L
\phi_{r,N,j_\ell}\ast\psi_N(x_{\ell})
\nonumber
\end{eqnarray}
with $ \x_L =(x_{\ell})_{\ell=1}^{L}\in \R^L $ and 
$$
\phi_{r,N,j}(x)=r^{-N}{N\choose j}(r+x)^{j} (r-x)^{N-j},
$$
where $f\ast g$ stands for the convolution of $f$ and $g$.
Then by a simple observation we have the following.

\begin{lem}
\label{l:42}
Let $ g \in C_0^{\infty} (\R^L) $ with support in $[-r,r]^L$. 
Then 
\begin{eqnarray}& 
\lim_{N\to\infty}\int_{\R^L}\{|g-g_N|^2 
+\D(g-g_N, g-g_N) \} d\x_L = 0.
\nonumber
\end{eqnarray}
\end{lem}
Let $ L=dk$ and $ \x_L=\x_{dk}=(x_1,\ldots,x_k)$, where $ x_i\in S$. 
If a function $ g (x_1,\ldots,x_k)$ is symmetric in 
$ (x_1,\ldots,x_k)$, then we can and do regard $ g $ as a function 
defined on the configuration space $ \mathcal{M} $ over $ S $ 
with support in $ \{ \xi (S) = k \} $. 
If $ g (x_1,\ldots,x_k)$ is symmetric in 
$ (x_1,\ldots,x_k) $, where $ x_j\in S $, then so is $ g_N$. 
Hence, we deduce from Lemma \ref{l:42} the following. 
\begin{lem}
\label{l:43}
Let $L=dk$. 
Let $ g \in C_0^{\infty} (\R^L) $ be symmetric in $ (x_1,\ldots,x_k) $
with support in $[-r,r]^L$. 
Let $ g $ and $ g_N $ be regarded as functions on $ \mathcal{M} $ 
with support in $ \{ \xi (S) = k \} $ as above. Then 
\begin{align}\label{:43a}&
\lim_{N\to\infty} \|g - g_N \|_1 = 0
.\end{align}
\end{lem}

\bigskip 

\noindent {\bf Proof of Theorem \ref{th:main}}. 
We prove  only 
$(\cE^{\mu}, \cP^{\mu})= (\cE^{\mu}, \cD^{\mu})$ 
because the proof of the rest is similar. 
Let 
$$ S_{\mathrm{off}}^n = 
\{ (x_1,x_2,\dots,x_n)\in S^n ; x_i\not=x_j \text{ for }i\not=j \} 
.$$
Let $f\in\cD_\infty$ such that $ f (\xi ) = f (\xi_K )$ 
with compact set $ K = [-r,r]^d $. 
For such an $f$ we introduce continuous functions 
$\widehat{f}_{n}$, $n\in\N_0$ such that $ \widehat{f}_{0} = \check{f}_0 $ 
and that, for $ n \in \mathbb{N}$ and 
$ (x_1,x_2,\dots,x_n) \in S_{\mathrm{off}}^n$, 
\begin{eqnarray}
\label{fhat}
&& \widehat{f}_{n}(x_1,x_2,\dots,x_n) 
\nonumber\\
&& = \sum_{k=0}^n (-1)^{n-k} \sum_{\{i_1,\dots,i_k\}\subset \{1,2,\dots,n\}}\check{f}_{k}(x_{i_1},\dots,x_{i_k}) 
\nonumber
.\end{eqnarray}
The values of $ \widehat{f}_{n}$ on $ \{S_{\mathrm{off}}^n\}^c$ 
are defined by continuity. 
Then $\widehat{f}_{n}$ is a smooth symmetric function on $S^n$ 
vanishing out of $K$ for $ n \ge 1 $. 
Note that $\widehat{f}_n $ is the M\"obius transformation of 
$\check{f}_k$, $k=0,1,\dots,n$. Then we easily deduce that, 
for $ (x_1,x_2,\dots,x_n) \in S_{\mathrm{off}}^n$, 
\begin{eqnarray} \label{:44b}
&& \ \check{f}_k(x_1,x_2,\dots,x_k) 
\\
= &&\sum_{n=0}^k \sum_{\{i_1,\dots,i_n\}\subset \{1,2,\dots,k\}} 
 \widehat{f}_{n}(x_{i_1},x_{i_1},\dots,x_{i_n})
\nonumber
.\end{eqnarray}
This implies that $\check{f}_k - \check{f}_0$ 
can be represented by a linear combination of 
symmetric smooth functions vanishing out of $ K^k $.

In \eqref{:44c} and \eqref{:44d} below, 
$ (x_{i_1},x_{i_2},\dots,x_{i_n})$ and 
$ (y_1,y_2,\dots,y_{n })$ are taken to be in $  S_{\mathrm{off}}^n$. 
The equalities can be exteded to $  \{S_{\mathrm{off}}^n\}^c$ 
by continuity of the functions. 

Let $ k= \xi_{K}(S) $ and write 
$ \xi_{K} = \sum_{i=1}^{k} \delta_{x_i}$. 
Then from \eqref{:44b} we can rewrite $f$ as 
\begin{align} \label{:44c}
 f (\xi) &= \check{f}_{k}(x_1,x_2,\dots,x_{k})
\\\notag 
&= 
\sum_{n =0}^{\xi_{K}(S)} 
\sum_{\eta \prec \xi_K }
\widehat{f}_{n }(y_1,y_2,\dots,y_{n }) 
.\end{align}
Here $n = \eta (S)  $, 
$ \eta = \sum_{i=1}^{n } \delta_{y_i}$, and 
$ \eta \prec \xi _{K} $ means that $ \eta (A) \le \xi _{K}(A) $
 for all $ A $. 
We note that the right-hand side can be regarded as a symmetric function of 
$ (x_1,\ldots,x_n)$ by construction. 
For $ m \in \N $ we put
\begin{align}\label{:44d}&
f_{[m]}(\xi)= 
\sum_{n =0}^{m} 
\sum_{\eta \prec \xi_K }
\widehat{f}_{n }(y_1,y_2,\dots,y_{n })
.\end{align}
Let $\varepsilon > 0$ be arbitrary. 
Then from Lemma \ref{l:41} and \eqref{:44d}  
we can take $\a\in {\bf A}$ and $a_r \le m\in \N$ such that 
\begin{align}\label{:qqq}
 \| f - f_{[m]} \|_1 
 \le \, & \|(1 - {\cal X}[\a]) ( f - f_{[m]} ) \|_1
<\, \varepsilon 
.\end{align}

From Lemma \ref{l:43}, 
we approximate the symmetric function $\widehat{f}_n$ in \eqref{:44d}  
by a polynomial $\widehat{F}_n$. 
Hence, for any $\varepsilon >0$, 
we can take polynomials $F_n$ such that
\begin{equation} \label{est:A1} 
\|f_{[m]}- \sum_{n=0}^{m} F_n \|_1 <\varepsilon 
.\end{equation}

Results  
\eqref{:qqq} and \eqref{est:A1} complete the proof. 
\qed

\vskip 3mm
\noindent{\bf Acknowledgement}

H.O. is supported in part by a Grant-in-Aid for Scientific Research (KIBAN-A, No. 24244010)
of the Japan Society for the Promotion of Science.
H.T. is supported in part by a Grant-in-Aid for Scientific Research (KIBAN-C, No. 23540122)
of the Japan Society for the Promotion of Science.

\end{document}